\documentclass[12pt]{amsart}
\usepackage[margin=1in]{geometry}

\usepackage{graphicx}
\usepackage{amsthm}
\usepackage{amssymb}
\usepackage{xcolor}

\usepackage{enumitem}
\usepackage{hyperref}

\newtheorem{theorem}{Theorem}
\newtheorem{prop}[theorem]{Proposition}

\theoremstyle{definition}

\theoremstyle{remark}

\numberwithin{theorem}{section}
\numberwithin{equation}{section}

\newcommand{\Z}{\mathbb{Z}}

\renewcommand{\phi}{\varphi}

\newcommand{\set}[1]{\left\{#1\right\}}

\begin{document}

\title{Addendum to ``Fuglede's conjecture fails in 4 dimensions over odd prime fields''}

\author{Samuel J.~Ferguson}
\address{Courant Institute of Mathematical Sciences\\ New York University\\ 251 Mercer St.\\ New York, NY 10012, USA}
\email{ferguson@cims.nyu.edu}

\author{Nat Sothanaphan}
\address{Courant Institute of Mathematical Sciences\\ New York University\\ 251 Mercer St.\\ New York, NY 10012, USA}
\email{natsothanaphan@gmail.com}

\begin{abstract}
We describe computer programs accompanying our paper \cite{Prev} and show that running them suffices to verify Fuglede's conjecture in $\Z_2^5$ and $\Z_2^6$.
\end{abstract}

\maketitle

\section{Overview}

Accompanying our paper \cite{Prev} are two computer programs, \verb|DephasedRank| and \verb|CheckSpecTile|.
These, together with theoretical results from \cite{Prev}, form a complete proof that Fuglede's conjecture holds in $\Z_2^5$ and $\Z_2^6$. The goal of this addendum is to show that running the programs suffices to verify the conjecture. Since implementation details are given in text files and comments in code files, we mainly give the big picture of the programs. Our code can be accessed at

\begin{center}
\url{https://github.com/natso26/FugledeZ_2-d}.
\end{center}

The code is written in Python 3, with MATLAB API for Python required for computing ranks in $\Z_2$ for \verb|DephasedRank|.

\section{DephasedRank}
\label{sec:dephased}

The goal of \verb|DephasedRank| is to prove the following proposition.

\begin{prop}
\label{prop:maindephased}
There is no spectral set in $\Z_2^5$ or $\Z_2^6$ whose size is not a power of $2$.
\end{prop}

\noindent
Since, clearly, there is no tile in $\Z_2^5$ or $\Z_2^6$ whose size is not a power of $2$, Proposition \ref{prop:maindephased} reduces Fuglede's conjecture on such groups to the equivalence of spectrality and being a tile for sets whose sizes are powers of $2$. This latter statement will be verified by \verb|CheckSpecTile|.

By \cite[Thm. 4.7]{REU}, Proposition \ref{prop:maindephased} is equivalent to the nonexistence of a log-Hadamard matrix with entries in $\Z_2$ and size $m \times m$, where $m$ is not a power of $2$, of rank at most $6$. Our \cite[Thm. 4.3]{Prev} says that a \emph{dephased} log-Hadamard matrix has the lowest rank among all log-Hadamard matrices in the same equivalence class. Hence, Proposition \ref{prop:maindephased} is equivalent to the following.

\begin{prop}
\label{prop:dephasedreduce}
The rank of \emph{some} dephased log-Hadamard matrix in every equivalence class of log-Hadamard matrices with entries in $\Z_2$ and size $m \times m$, where $m < 64$ is not a power of $2$, is more than $6$.
\end{prop}

Note the importance of ``some'' in Proposition \ref{prop:dephasedreduce}: it implies that computing the rank of any dephased matrix in a given equivalence class suffices.
Further restrictions on $m$ in Proposition \ref{prop:dephasedreduce} are possible. By \cite[Thm. 1.1(b),(c)]{REU}, we can assume that $m$ is a multiple of $4$ and is less than $32$. Moreover, since Fuglede's conjecture holds in $\Z_2^4$ \cite[Prop. 5.2]{Prev}, it is enough to check when $m > 16$. Thus, we arrive at the following reduction of Proposition \ref{prop:maindephased}.

\begin{prop}
\label{prop:dephasedreduce2}
The rank of \emph{some} dephased log-Hadamard matrix in every equivalence class of log-Hadamard matrices with entries in $\Z_2$ and size $m \times m$, where $m=20,24,28$, is more than $6$.
\end{prop}

To prove Proposition \ref{prop:dephasedreduce2}, we employ information in Sloane's online library of Hadamard matrices \cite{Slo}, which lists a Hadamard matrix from every equivalence class up to size $28 \times 28$. It is fortunate that this barely covers every multiple of $4$ up to less than $32=2^6/2$, which allows us to handle $\Z_2^6$. Indeed, the numbers of equivalence classes dramatically grow beyond this point: there are $487$ equivalence classes of Hadamard matrices of order $28$, but over $10$ million of order $32$. So handling $\Z_2^7$ is likely impossible without new theoretical results.

There are 3, 60 and 487 equivalence classes of Hadamard matrices of orders 20, 24 and 28, respectively. For each equivalence class, Sloane's library gives a Hadamard matrix in that class. \verb|DephasedRank| can read such a matrix (fetch the file directly from the URL); convert it to the corresponding log-Hadamard matrix by changing $1$ and $-1$ to $0$ and $1$, respectively; dephase it; and compute its rank in $\Z_2$. The program can also read many files at once, which allows us to process the 60 matrices of order 24 and 487 matrices of order 28 simultaneously. The rest are implementation details, which are provided in the text and code files on GitHub.

The results are as follows. For the $3$ equivalence classes of matrices of order 20, all dephased matrices have rank 18. For the 60 classes of order 24, all have rank 11; and for the 487 classes of order 28, all have rank 26. Thus, the ranks are all greater than 6, which proves Proposition \ref{prop:dephasedreduce2} and hence Proposition \ref{prop:maindephased}.

It is curious that dephased matrices in all equivalence classes of log-Hadamard matrices of the same order seem to have the same rank. Indeed, this is one of our conjectures in \cite[Sec. 5]{Prev}, which states that this is true as long as the order is not a power of 2. (For order 16, for example, the ranks can be either $4$, $5$, $6$ or $7$.) If this conjecture holds, then \verb|DephasedRank| would also be able to handle $\Z_2^7$, since Sloane's library provides at least one matrix of every order up to order $256 > 2^7/2$.

\section{CheckSpecTile}

We now move on to \verb|CheckSpecTile|, which does more heavyweight work. The goal of this program is to prove the following proposition.

\begin{prop}
\label{prop:maincheckspec}
A set in $\Z_2^5$ or $\Z_2^6$ whose size is a power of $2$ is spectral if and only if it is a tile.
\end{prop}

First, some theoretical reductions are possible. By \cite[Thm. 1.1(e)]{REU}, it suffices to consider sets of sizes $4$, $8$ in $\Z_2^5$ and sizes $4$, $8$ and $16$ in $\Z_2^6$. The next reduction is that subsets of $\Z_2^5$ need not be considered at all, according to the proposition below.

\begin{prop}
\label{prop:affinereduce}
Let $m<n$. A subset of $\Z_2^m$ is spectral or is a tile if and only if it is spectral or is a tile, respectively, when thought of as a subset of $\Z_2^n$ by identifying $\Z_2^m$ with $\Z_2^m \times \set{0}^{n-m} \subseteq \Z_2^n$.
\end{prop}

The spectrality part of Proposition \ref{prop:affinereduce} is a special case of the affine restriction proposition in our working paper \cite[Cor. 3.16 or 3.17]{Riesz}. The tiling part is straightforward.

Next, in the case of subsets of size 4 of $\Z_2^6$, such a subset must lie in an affine space of dimension $3$, so affine invariance \cite[Cor. 3.2(1) and 4.3(b),(c)]{REU} together with Proposition \ref{prop:affinereduce} reduce our question to Fuglede's conjecture in $\Z_2^3$, which is true. Therefore, we are left with the following proposition to prove.

\begin{prop}
\label{prop:newcheckspec}
A set in $\Z_2^6$ of size $8$ or $16$ is spectral if and only if it is a tile.
\end{prop}

Proving Proposition \ref{prop:newcheckspec} is a huge amount of work: there are $\binom{64}{16} > 10^{14}$ subsets of size $16$ in $\Z_2^6$, and for each set, the same number of candidates for its spectrum or tiling partner. So effort is required to make the computation manageable. We will first outline our verification approach, which results in an algorithm that is estimated to run for more than a decade. Then we describe our subsequent modifications, which speed the program up so that it can finish within a day.

\subsection{Approach}

We will focus on subsets of size $16$, since subsets of size $8$ are easier. The first idea is to again use affine invariance to assume that the first $7$ points are
$$(0,0,0,0,0,0), \,\, (0,0,0,0,0,1), \,\, (0,0,0,0,1,0), \,\, \dots, \,\, (1,0,0,0,0,0).$$
We can do this as long as the set does not lie in any proper affine subspace. This assumption is justified because if the set does lie in a proper affine subspace, then affine invariance and Proposition \ref{prop:affinereduce} reduce the question to Fuglede's conjecture in $\Z_2^5$ for subsets of size $16$, which is true by \cite[Thm. 1.1(e)]{REU}.

For subsets of size $8$, a similar argument yields that we can assume the first $7$ points as above, or we must consider subsets of size $8$ in $\Z_2^5$ with the first $6$ points assumed in a similar manner. There are $\binom{57}{1}=57$ and $\binom{26}{2}=325$ subsets in each case, which are rather small numbers. The case of subsets of size $16$ leaves us with $\binom{57}{9} > 10^9$ subsets, a number whose size will cause us some difficulty.

We now explain how to determine whether a given set is spectral or is a tile. Clearly, enumerating the $\binom{64}{16}$ potential spectra or tiling partners is a not a good idea. Instead, we use an algorithm based on \emph{clique finding}, a famous problem in computer science. This clique approach is also taken by Siripuram et al. \cite{Sir} in studying spectrality in groups $\Z_N$.

We first illustrate the clique approach for the tiling part of the problem. To see whether a set tiles, put out all $64$ translates of the set as nodes in a graph, and draw an edge between two nodes if the two copies do not overlap. Then the set tiles if and only if the constructed graph contains a clique of size $64/16=4$. Moreover, we can always assume that the original copy is part of the clique, so it suffices to find a $3$-clique among nodes connected to the original copy.

The spectrality part is a little trickier. If $E$ is our set and $B$ is its spectrum, then the exponentials corresponding to elements of $B$ must be orthogonal on $E$. So we put the $64$ exponentials as nodes in a graph, and draw an edge between two nodes if the two exponentials are orthogonal on $E$. The question then reduces to finding out whether a $16$-clique exists in such a graph. By translation invariance \cite[Cor. 4.3(c)]{REU}, the exponential $1$ can be assumed to be in the clique, so it suffices to find a $15$-clique among its adjacent nodes.

Even though the clique problem is NP-complete, meaning there is unlikely to be an efficient algorithm for it, this does not cause much trouble for us in practice. Since a $3$-clique is small compared to the graph, and most graphs do not have such a large clique as a $15$-clique, a reasonable algorithm that takes these facts into account can perform rather well. (The clique finding part is not the major bottleneck of the program.)

More implementation details can be found in the folder \verb|CheckSpecTile/Original| on GitHub.

\subsection{Optimization}
\label{subsec:optimize}

Regardless, the approach taken above is still too slow for our purpose. With 59 ms for the computation for each subset, the large number of subsets means that the entire enterprise will take about 17 years. So we perform two kinds of optimizations. First, we reduce the computation time for each subset to $170$ $\mu$s, an improvement by a factor of $360$. Then, we use symmetry considerations to reduce the number of subsets by a factor of $30$. The resulting code is $4$ orders of magnitude faster than the original one and takes only $13$ hours to run.

The first kind of optimizations is done through code profiling to identify bottlenecks. We find the computation of the $64 \times 64$ adjacency matrix for the graph to be a bottleneck. To remedy this, the matrix is never explicitly computed, but translation invariance is used so that only nodes adjacent to node $0$ are computed. More precisely, node $x$ is adjacent to node $y$ if and only if node $0$ is adjacent to node $x-y$, the subtraction being done componentwise modulo $2$. The computation of which nodes are adjacent to node 0 for both spectrality and tiling is still a bottleneck, so the code for these parts is extensively rewritten. Finally, the clique finding algorithm is revised so as to reduce the amount of recursion being done and run faster in some parts.

The authors considered using an external library to solve the clique problem at some point. However, although it may solve the clique problem faster, the generation of the graph is itself a bottleneck, so this is ultimately slower than our approach which keeps only nodes adjacent to node 0. The authors also considered going from a Python to a C (a faster language) implementation, but this turns out to be unnecessary.

More details can be found in code files in the folder \verb|CheckSpecTile/Fast| on GitHub.

The second kind of optimizations is concerned with coordinate symmetry: permuting coordinates preserves the properties of being spectral and being a tile. Apart from the $7$ fixed points, let the remaining $9$ points be $x_1 < x_2 < \dots < x_9$, arranged in the lexicographical order. Then we enumerate all $9$-tuples $(x_1,x_2,\dots,x_9)$ in the lexicographical order. The rule is that if a $9$-tuple can be coordinate-permuted into another tuple that comes before it in our order, then it can be safely ignored. The reason is that some coordinate-permuted version of the tuple would have already been examined.

We do not perform a perfect check of the rule, as it will take too much computation time, but some heuristics are implemented to throw out some $x_1$; some $x_2$ given $x_1$; some $x_3$ given $x_1$ and $x_2$; and some $x_i$, $i \geq 4$, given $x_1$, $x_2$ and $x_3$. There are $1$ heuristic concerning $x_1$, $3$ concerning $x_2$, $4$ concerning $x_3$, and $4$ concerning $x_i$ when $i \geq 4$. Some examples of heuristics are as follows.
\begin{itemize}
    \item All $0$ entries of $x_1$ must come before all $1$ entries. For example, $x_1 = (0,0,1,0,1,1)$ can be thrown out.
    \item Restricting to coordinates of $x_1$ that are 0, the coordinates of $x_2$ must be nondecreasing. For example, $x_1 = (0,0,0,1,1,1)$ and $x_2 = (1,0,1,0,1,1)$ are not possible.
\end{itemize}
A list of heuristics is provided in the Appendix, and a full description and proof are in the file \verb|CheckSpecTile/Fast/generate.py|. 

\section{Discussion}

Together, the programs \verb|DephasedRank| and \verb|CheckSpecTile| settle Fuglede's conjecture in $\Z_2^5$ and $\Z_2^6$. A natural question is whether the same approach can be used to tackle the case of $\Z_2^7$. We find this to be unlikely without further theoretical results. For \verb|DephasedRank|, a proof of the conjecture discussed in the last paragraph of Section \ref{sec:dephased} is desirable and will indeed allow the program to work in the case of $\Z_2^7$.

For \verb|CheckSpecTile|, however, the case of $\Z_2^7$ seems out of reach. With $8$ fixed points, the number of subsets of size $32$ is $\binom{120}{24} > 10^{25}$, which is enormous. Coordinate symmetry, which can reduce this number by at most a factor of $7!=5040$, would be far from enough. Moreover, clique finding in this case becomes more complicated. The computation for each subset takes $520$ $\mu$s, which is about $3$ times that of the case of $\Z_2^6$.

It should be remarked that the NP-completeness of the clique problem means that the computation time for each subset may dramatically increase as we go on to $\Z_2^8$ and $\Z_2^9$. Relatedly, Kolountzakis and Matolsci \cite[Sec. 4]{KM} also show that a problem related to spectrality and tiling is NP-complete.

We conjecture that, ultimately, Fuglede's conjecture holds in $\Z_2^d$ for all  $d \leq 9$. We already know that it fails in $\Z_2^{10}$ \cite[Thm. 1.2]{Prev}, so we think that this result is sharp. The reason is that the $\Z_2^{10}$ counterexample is a log-Hadamard matrix of size $12 \times 12$ with rank $10$. However, as the matrix size increases, the rank seems to also increase in general and seems unlikely to fall below $10$ again (when the order is not a power of two). This gives some heuristics as to why $10$ may be the smallest dimension such that Fuglede's conjecture fails in $\Z_2^d$.

\section*{Appendix}

We list here the heuristics used in Section \ref{subsec:optimize} to reduce the number of subsets to check for tiling and spectral properties.

For $x_1$:
\begin{itemize}
\item $x_1$ must be of the form $(0,\dots,0,1\dots,1)$, i.e. the coordinate values are nondecreasing.
\end{itemize}

For $x_2$:
\begin{itemize}
\item $x_2 > x_1$.
\item The number of coordinates of $x_2$ that are $1$ must be at least the number of coordinates of $x_1$ that are $1$.
\item Restricting to the coordinates for which $x_1$ is $0$, the coordinates values of $x_2$ are nondecreasing. Likewise, restricting to the coordinates for which $x_1$ is 1, the coordinate values of $x_2$ are nondecreasing.
\end{itemize}

For $x_3$:
\begin{itemize}
\item $x_3 > x_2$.
\item The number of coordinates of $x_3$ that are $1$ must be at least the number of coordinates of $x_1$ that are $1$.
\item Restricting to the coordinates for which $x_1$ is $0$, the number of coordinates of $x_3$ that are 1 must be at least the number of coordinates of $x_2$ that are 1.
\item For two consecutive coordinates such that the values of $x_1$ on these coordinates are the same, and likewise for $x_2$, the coordinate values of $x_3$ are nondecreasing. In other words, we cannot have  $x_1 = (\dots,a,a,\dots)$, $x_2 = (\dots,b,b,\dots)$, and $x_3 = (\dots,1,0,\dots)$.
\end{itemize}

For $x_i$ when $i \geq 4$:
\begin{itemize}
\item $x_i > x_3$.
\item The number of coordinates of $x_i$ that are $1$ must be at least the number of coordinates of $x_1$ that are $1$.
\item Restricting to the coordinates for which $x_1$ is $0$, the number of coordinates of $x_i$ that are $1$ must be at least the number of coordinates of $x_2$ that are $1$.
\item Let the first $k$ coordinates of $x_2$ be 0. Then restricting to the first $k$ coordinates, the number of coordinates of $x_i$ that are $1$ must be at least the number of coordinates of $x_3$ that are $1$.
\end{itemize}
\end{document}